\documentclass[final,onefignum,onetabnum]{siamart171218}
\usepackage{amsmath}
\usepackage{amsfonts}
\usepackage{algorithmic}
\usepackage[caption=false]{subfig}
\usepackage[colorinlistoftodos, textwidth=30mm]{todonotes}
\presetkeys{todonotes}{backgroundcolor=white,linecolor=black}{}
\usepackage[normalem]{ulem}

\usepackage{color}

\newcommand{\added}[1]{#1}
\newcommand{\ssout}[1]{}
\newcommand{\changed}[1]{#1}
\newcommand{\removed}[1]{}
\newcommand{\done}[2][]{}
\newcommand{\mdone}[1]{}

\newcommand{\R}{\mathbb{R}} 
\newcommand{\cP}{ {\bf P} } 

\newcommand{\cV}{{\mathcal{V}}} 

\newcommand{\Le}{{\bf L}} 
\newcommand{\cW}{{\mathcal{W}}}
\newcommand{\bcW}{\pmb{\mathbf{\cW}}}

\newcommand{\ds}{{\displaystyle}}

\newcommand{\D}{{\bf D}}  
\newcommand{\bstar}[1]{ {\bf {#1}^*}}

\def \ds          {\displaystyle}

\def \cT          {{\cal T}}

\def \bg          {{\bf g}}

\renewcommand{\theequation}{\thesection.\arabic{equation}}

\newcommand{\WDNote}[1]           
{\textcolor{red}{#1}\marginpar{\textcolor{red}{WD $\longleftarrow$}}}
\newcommand{\SNote}[1]           
{\textcolor{blue}{#1}\marginpar{\textcolor{blue}{SW $\longleftarrow$}}}
\newcommand{\Fig}[1]{Fig.~\ref{#1}} 
\newcommand{\Figure}[1]{Figure~\ref{#1}} 
\newcommand{\Sec}[1]{Sec.~\ref{#1}} 

\graphicspath{{./figures/}}

\DeclareMathOperator*{\supp}{supp}

\title{\bf Simultaneous Sensing Error Recovery and Tomographic Inversion
Using an Optimization-based Approach }

\author{Anthony P. Austin\thanks{Mathematics and Computer Science Division,
Argonne National Laboratory.\
({\tt austina@anl.gov}, {\tt wendydi@anl.gov}, {\tt leyffer@anl.gov}, {\tt 
wild@anl.gov}.)}
  \and Zichao (Wendy) Di\footnotemark[3]
  \and Sven Leyffer\footnotemark[3]
\and Stefan M.~Wild\footnotemark[3]
}

\begin{document}

\maketitle
\done[inline]{\em We have used the following convention to highlight our 
changes in response to the referees' comments: red text indicates modification, 
purple text indicates addition, and margin parameters indicate responses to 
specific comments.}

\begin{abstract}
Tomography can be used to reveal internal properties of a 3D object using any penetrating wave. Advanced tomographic imaging techniques,
however, are vulnerable to both systematic and random errors associated with
the experimental conditions, which are often beyond the capabilities of the state-of-the-art reconstruction techniques such as regularizations.  Because they can lead to reduced spatial
resolution and even misinterpretation of the underlying sample structures,
these errors present a fundamental obstacle to full realization of the
capabilities of next-generation physical imaging.
In this work, we develop efficient and explicit recovery schemes of the most
common experimental error:  movement of the center of rotation during
the experiment.  We formulate new physical models to capture the experimental setup,
and we devise new mathematical optimization formulations for reliable inversion of
complex samples. We demonstrate and validate the efficacy of our approach on
synthetic data under known perturbations of the center of rotation.
\end{abstract}


\section{Introduction}
\label{sec:intro}

Tomographic imaging has had a revolutionary impact on medicine, physics, and
chemistry.  Even so, the problem of reconstructing an image from tomographic
data remains challenging in many interesting cases, such as when the amount of 
available data is
limited and/or the problem is ill posed (in the sense that canonical metrics 
used to assess the discrepancy between a reconstruction and the
measured data generally possess many local minima).  This ill posedness makes
reconstructions susceptible to experimental errors, in particular, to errors
stemming from mismatches between the experimental configuration and the
assumptions of the measurement process. Recovering such errors is crucial for 
realizing \changed{the gains from improvements in} measurement and experimental hardware, such as 
the improved
resolution promised by brighter, more coherent next-generation light sources.
Sample drift \cite{mlodzianoski2011sample} and beam drift
\cite{tripathi2014ptychographic} are two fundamental sources of error. These errors can often arise from the drift of the center of rotation (CoR) of the imaging stage
during data acquisition \cite{azevedo1990calculation}, which is the focus of this work.

In computerized tomography, the object or sample being imaged is placed on a stage and
irradiated with parallel beams of x-rays.  As the rays pass through the object,
they are partially absorbed according to the object's composition.  Radiation
that passes through the object unabsorbed is collected by a detector, producing
an x-ray ``shadow'' or ``projection'' of the object.  By rotating the object
and repeating the process, one obtains projections of the object as illuminated
from several different directions. One then uses these projections to
reconstruct an image of the object and its interior \cite{kalender2006x}.  

\begin{figure}[t]
  \centering
  \includegraphics[scale=0.8]{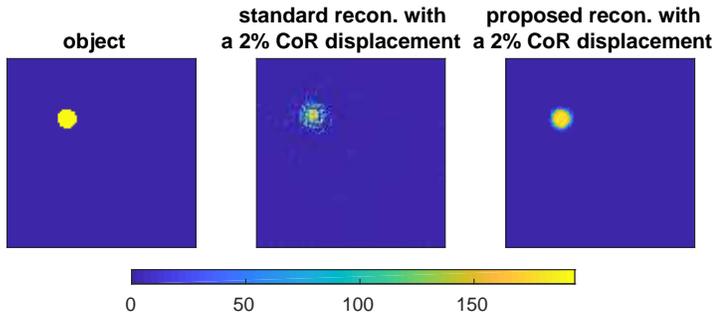}
  \caption{Consequences of failing to account for CoR shift during tomographic
reconstruction.  In this experiment, the CoR was displaced by 2\% with respect 
to the size of the object domain, 
\changed{and} 
the standard reconstruction without CoR recovery 
(middle) looks very different from the
original \added{spherical} 
\changed{object} 
(left). Our proposed reconstruction with CoR recovery is shown on the right, which almost perfectly recovers the object.}
  \label{fig:corError}
\end{figure}

Imperfections in the experimental apparatus can cause the CoR of the imaging
stage to vary slightly over the course of this process.  When this happens, the
projections from different angles will not be properly aligned relative to one
another, and failing to account for this misalignment during reconstruction may result in a
smeared image (loss of resolution) or, worse, an image that is completely
incorrect.  Moreover, even a small error in the CoR can yield a large error in
the reconstruction, as is illustrated by the experiment of
\Fig{fig:corError}.  The image on the left shows the object---a small,
solid circle---and the image \changed{in the middle} shows the result of a standard tomographic
reconstruction \added {without addressing the CoR shift} 
from a sequence of measurements in which
the CoR is perturbed by a 2\% displacement relative to the size of 
the object domain.  The reconstructed image obtained by not correcting for this 
perturbation looks very different from the true object. \added{We also demonstrate the reconstruction of the approach we propose in this paper on the right side, which almost perfectly resembles the true object.}

In the past, it has been possible largely to ignore errors like these, as the
drifts have been small compared with the widths of the beams used to illuminate
the object. With demand for increasingly finer resolutions leading \changed{from}
higher-quality light sources with narrower beams, however, these errors can no 
longer be
ignored. Potential approaches to this problem include the incorporation
of prior knowledge into the reconstruction process using Bayesian frameworks
\cite{hanson1983bayesian,sauer1994bayesian,spantini2015optimal} and the use of
regularizers to promote sparsity or smoothness of the reconstruction
\cite{hamalainen2013sparse,girard1987optimal,lassas2004can,niinimaki2016multiresolution, hsieh2013recent}.
One can also collect additional data (beyond x-ray projections) that is less
sensitive to experimental errors and use this data to help with the reconstruction
\cite{di2016optimization}.  Unfortunately, the improvement in tomographic hardware 
has been so great that these generic approaches, which apply more
widely to problems other than CoR drift recovery, are not likely to be
sufficient. In the case of CoR drift, we show (in \Sec{sec:numerical}) that 
standard regularization approaches fail even when the regularization parameter 
value is chosen in an ideal way.

Much research effort has been 
\changed{devoted to} 
tackling the CoR 
\changed{recovery} 
problem. For example, \mdone{R1C1}\added{a typical approach 
to recovering a single, unknown CoR is to use a pair of projections that are 
reflections of one another to estimate a detector offset that can be used to 
shift the projections into the correct positions 
\cite{zitova2003image,min2012new}.  This approach is highly sensitive to the 
accuracy of the mirrored alignment, and its use of such a small number of 
projections to perform the correction makes it susceptible to noise and other 
effects of limited data, potentially resulting in a low-contrast recovery.  
Moreover, mirrored pairs of projections are not always available in practice.} 
Azvedo et al.\ \cite{azevedo1990calculation} proposed a method to estimate the 
CoR based on the preservation of photon counts passing across the sample, which 
requires almost perfect measurement. The most common technique is to compute
cross correlations between projections acquired from successive rotations
\cite{guckenberger1982determination,amat2010alignment,hayashida2010automatic};
however, this approach is limited to simple and relatively homogeneous samples,
since there is no good way to rank one feature as ``more important'' than another
when two projections contain very distinguishable features.  A similar approach
involves manual alignment of the projections using a known ``hotspot'' in the
object as a reference.  This process can be labor-intensive and cannot be used
if a hotspot cannot be identified.  Recently, a new technique known as
iterative reprojection has been introduced in which one recovers projection
alignments and reconstructs the image simultaneously using an iterative
procedure
\cite{tomonaga2014alternative,houben2011refinement,winkler2006accurate,gursoy2017rapid}.
The basic idea is to alternate between a few iterations of projection alignment
and a few iterations of reconstruction until a ``forward model'' of the
experimental setup and an ``inverse model'' of the reconstruction process are
consistent with one another.  This approach can yield better accuracy,
especially for noisy, limited data.

One deficiency of the iterative reprojection methods that have emerged thus far
is that the update of the projection alignment and the update of the
reconstruction are not fully coupled.  As a result, some
mismatch remains between the alignment of the projections and the experimental
configuration assumed by the reconstruction process.  In this work, we propose
to address this problem using a novel joint inversion framework based on
optimization in which we explicitly model and recover for CoR drift. \mdone{R1C1}\added{Our proposed model is flexible, making no assumptions about when and where the drift happens.}  Our new
approach is easier to automate than existing approaches and can be more robust
in the context of poor data quality and limited prior knowledge of the object
being imaged. 

In \Sec{sec:model}, we describe our mathematical forward models of the
CoR drift error in the experiment and show how to embed these models into the
reconstruction scheme.  In \Sec{sec:alg}, we describe our simultaneous
reconstruction approach for recovering the object and the experimental error,
including the formulation of the objective function.  In \Sec{sec:opt},
we describe the algorithm for solving the resulting optimization problem and
its complexity.  In \Sec{sec:numerical}, we present some numerical
illustrations comparing the performance of our simultaneous inversion method
with that of existing approaches using a few synthetic examples.  In
\Sec{sec:conc}, we summarize the proposed method and discuss a few
directions for future research.

For simplicity, we confine ourselves
to reconstructing 2D images, although our methods work just as well in 3D.

\section{Mathematical Model}
\label{sec:model}


In this section, we describe our model for the tomographic imaging process.
For further details on tomography, we refer the reader to  
\cite{kak_principles_1988}.

\subsection{Radon Transform}

The fundamental mathematical tool in tomography is the Radon transform
\cite{radon_determination_1986}, defined for a compactly supported function 
$f: \R^2\mapsto \R$ by
\begin{equation}\label{EQN:RT}
Rf(\tau, \theta) = \int_{-\infty}^\infty \int_{-\infty}^\infty f(x, y) \delta(\tau - x\cos \theta - y \sin \theta) \: dx \: dy,
\end{equation}
where $\delta$ is the Dirac delta function.  Throughout this article, we 
consider the restricted domain given by $\tau \in[0, \infty)$ and $\theta \in 
[0, 2\pi)$,
and we assume, as will always be the case in practice, 
that $f$ is
well-enough behaved that the integral makes sense. 
The Radon transform of a
function is frequently called its \emph{sinogram}.  
\begin{figure}[t]
  \centering
  \includegraphics[width=.7\linewidth]{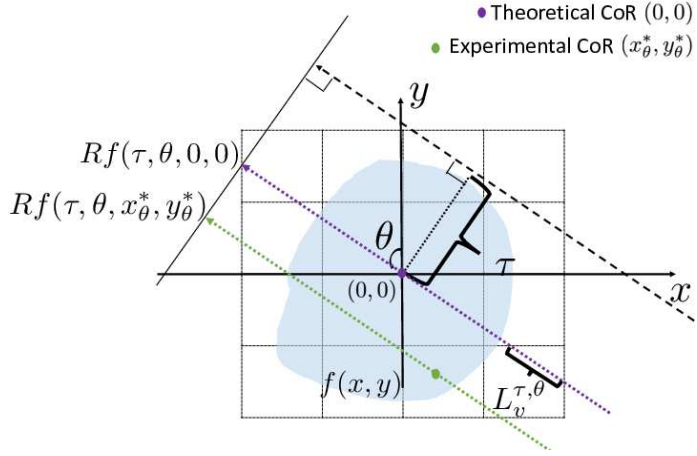}%
  \caption{Geometric sketch of the Radon transform, which maps $f$ from $(x,y)$ 
space to $(\theta,\tau)$ space. The purple line and the green line denote rotations of their previous position,
as $y$-axis, with respect to different CoR, respectively.}
  \label{fig:contRadon}
\end{figure}


The projections that are measured in the tomographic x-ray imaging process are 
values of the
Radon transform of the object's \emph{attenuation coefficient} $f(x,y)$, a function that
describes the propensity of the object to absorb x-rays at each point in the 
object's support.
The angle $\theta$ encodes the direction from which the beams approach the
object; and given that direction, $\tau$ encodes the position of the beam.  If
the CoR shifts from the origin to $(x_\theta^*, y_\theta^*) \in \R^2$, as the imaging stage rotates
through the angle $\theta$, the projections we measure come not from
\eqref{EQN:RT} but from the CoR-shifted Radon transform
\begin{equation}\label{EQN:ShiftedRT}
Rf(\tau, \theta, x_\theta^*, y_\theta^*) = \int_{-\infty}^\infty \int_{-\infty}^\infty f(x, y) \delta\bigl(\tau - x_\theta^* - (x - x_\theta^*) \cos \theta - (y - y_\theta^*) \sin \theta\bigr) \: dx \: dy.
\end{equation}
For a derivation of this equation, see Appendix \ref{SEC:ShiftedRTDerivation}.
As illustrated in \Fig{fig:contRadon}, $Rf(\tau, \theta,0,0)$ is the 
integral of $f$ along the purple line perpendicular
to the direction determined by the angle $\theta$ and at a distance $\tau$ from
the origin. Alternatively, the purple line is obtained by rotating from its previous position, which overlays with the $y$-axis, $\theta$ degrees with respect to CoR $(0,0)$. However, if the CoR is shifted from $(0,0)$, denoted by the purple dot, to the green dot, we obtain a different projection as $Rf(\tau,\theta,x_\theta^*,y_\theta^*)$.  

To recover for the difference in CoR, we need a way to convert
\eqref{EQN:ShiftedRT} back into \eqref{EQN:RT}.  This is easy: we have
\begin{equation}\label{EQN:RTFixCoR1}
Rf(\tau, \theta, 0, 0) = Rf(\tau - P_\theta, \theta, x_\theta^*, y_\theta^*),
\end{equation}
where
\begin{equation}\label{EQN:RTFixCoR2}
P_\theta = x_\theta^*(1 - \cos \theta) + y_\theta^* \sin \theta.
\end{equation}
The simple relationship \eqref{EQN:RTFixCoR1}--\eqref{EQN:RTFixCoR2} forms the
basis for our algorithm.  We rotate the imaging stage through an angle $\theta$
with the intent of measuring $Rf(\tau, \theta, 0, 0)$, but we instead measure
$Rf(\tau, \theta, x_\theta^*, y_\theta^*)$ due to drift of the CoR.  The identity
\eqref{EQN:RTFixCoR1} shows that all we need to do is translate (in $\tau$) the
Radon transform that we measured by an amount $P_\theta$, which is related to
the true CoR by \eqref{EQN:RTFixCoR2}.  As we later show, given 
sufficient data, we can estimate the drift-induced 
$P_\theta$ using an optimization procedure.

\subsection{Discrete Tomography}

In practice, we cannot recover the desired object property (e.g., attenuation coefficient) at all points in
space.  Instead, we discretize the space (containing the compact object) being 
imaged into \changed{$N\times N$ pixels and vectorize it to an array $\cV$}.  Let $\cW_v$ be the value of the object property we intend to recover on the pixel
$v \in \cV$, and let $\bcW = \{\cW_v \: : \: v \in \cV\}$ denote the
discretized image.  Given parameters $\tau$ and $\theta$, we calculate a
discrete Radon transform of $\bcW$ via
\[
R_{\tau, \theta}(\bcW) = \sum_{v \in \cV} L_v^{\tau, \theta} \cW_v,
\]
where $L_v^{\tau, \theta}$ is the length of the intersection of the beam
described by $\tau$ and $\theta$ with the pixel $v$, see 
\Fig{fig:contRadon}.  
It is this relationship that we must invert to reconstruct our image:
knowing $R_{\tau, \theta}(\bcW)$, we wish to find $\cW_v$ for each pixel $v$.

In addition to needing to consider the image discretely, we must also
``discretize the beams'':  we have access only to $R_{\tau, \theta}$ for a
limited number of values of $\tau$ and $\theta$.  The values of $\tau$ we use
are fixed and equally spaced and are identical for each $\theta$.  We denote
the collection of values of $\tau$ (sometimes called ``beamlets'') by $\cT$ and
the collection of values of $\theta$ by $\Theta$; \added{accordingly, $N_\tau = |\cT|$ is the number of beamlets, and $N_\theta = |\Theta|$ is the number of angles.} In reality, the resolution of $\cT$ is decided by the energy level of the radiation source and by  the detector resolution \cite{bonse1996x}.

Thus, for each angle $\theta \in \Theta$, our measurement apparatus ideally
produces a set of samples of $Rf(\tau, \theta)$, equally spaced in $\tau$,
which we take as values of $R_{\tau, \theta}(\bcW)$.  Because of the drift in the
CoR, however, our samples are actually from $Rf(\tau, \theta, x_\theta^*, y_\theta^*)$
instead.  We can account for this drift by recomputing the
$L_v^{\tau, \theta}$ to be consistent with the change in CoR, but we wish to
avoid doing so because calculating these values (potentially many times) is 
expensive.  Instead, we calculate
them once, assuming that the CoR is at the origin, and recover for the change
in CoR by translating the projections.

\subsection{Aligning the Projections}
\label{SSEC:Aligning}

In principle, translating $Rf(\tau, \theta, x_\theta^*, y_\theta^*)$ in the $\tau$ variable
by $P_\theta$ is an easy task---the answer is just $Rf(\tau - P_\theta, \theta,
x_\theta^*, y_\theta^*)$---but since we have access only to samples of $Rf(\tau, \theta, x_\theta^*,
y_\theta^*)$ at $\tau \in \cT$, we cannot do this.  Since the samples come from
equally spaced points, a natural idea is to effect the translation by using the
discrete Fourier transform, which is mathematically equivalent to forming a
trigonometric interpolant to $Rf(\tau, \theta, x_\theta^*, y_\theta^*)$ through the points
$\tau \in \cT$ and translating the interpolant.  This approach works, but 
the images recovered in this way in practice are contaminated with ringing
artifacts \cite[pp.~209]{bracewell1986fourier}.  To fix this problem, we apply a low-pass filter---a normalized Gaussian
filter with standard deviation $\sigma$---to damp the high-order Fourier
coefficients.

An alternative way to understand the translation process is as follows.
Translating $Rf(\tau, \theta, x_\theta^*, y_\theta^*)$, viewed as a function of $\tau$ only,
by $P_\theta$ is equivalent to convolving $Rf$ with a Dirac delta function
centered at $P_\theta$:  $Rf(\tau - P_\theta) = (Rf \ast
\delta_{P_\theta})(x)$, where $\delta_{P_\theta}(\tau) = \delta(\tau -
P_\theta)$.  We regularize $\delta_{P_\theta}$ by replacing it with a Gaussian,
\begin{equation}\label{EQN:Gaussian}
\delta_{P_\theta, \sigma}(\tau) = \frac{1}{\sigma\sqrt{2\pi}} e^{-\frac{(\tau - P_\theta)^2}{2\sigma^2}},
\end{equation}
and instead compute the convolution $(Rf \ast \delta_{P_\theta,
\sigma})(\tau)$, which we do by sampling $Rf$ and $\delta_{P_\theta, \sigma}$
on the same equally spaced grid and using the discrete Fourier transform.

How should one choose the hyperparameter $\sigma>0$?  Let $\widetilde{Rf}(\tau, 
\theta, 0, 0)$ be our
approximation to $(Rf \ast \delta_{P_\theta, \sigma})(\tau)$%
\footnote{As in the preceding paragraph, the $Rf$ in $Rf \ast \delta_{P_\theta,
\sigma}$ in this definition is understood to be the CoR-shifted Radon
transform.  Our notation $\widetilde{Rf}(\tau, \theta, 0, 0)$ emphasizes our
hope that $\widetilde{Rf}(\tau, \theta, 0, 0) \approx Rf(\tau, \theta, 0, 0)$.}
obtained from the discrete Fourier transform.  By the error analysis presented
in Appendix~\ref{SEC:ErrorAnalysis}, if the (attenuation coefficient of the) 
object being imaged is twice continuously differentiable, then
\begin{equation}\label{EQN:TranslationError}
\widetilde{Rf}(\tau, \theta, 0, 0) = Rf(\tau, \theta, 0, 0) + O(\sigma^2) + O\left(\frac{1}{\sigma N_\tau}\right).
\end{equation}
\ssout{where $N_\tau = |\cT|$ is the number of beamlets.}  The second term on the
right-hand side of \eqref{EQN:TranslationError} represents the error incurred in 
our Gaussian
regularization of the delta function, while the third term  represents the
discretization error due to the limited number of 
\changed{measurements}.
A
trade-off between the two exists:  a narrower Gaussian (smaller $\sigma$) implies a
more faithful translation, but it requires more 
\changed{measurements} 
(greater $N_\tau$) to
approximate
accurately.  

\begin{figure}[t]
  \centering
  \includegraphics[width=.7\linewidth]{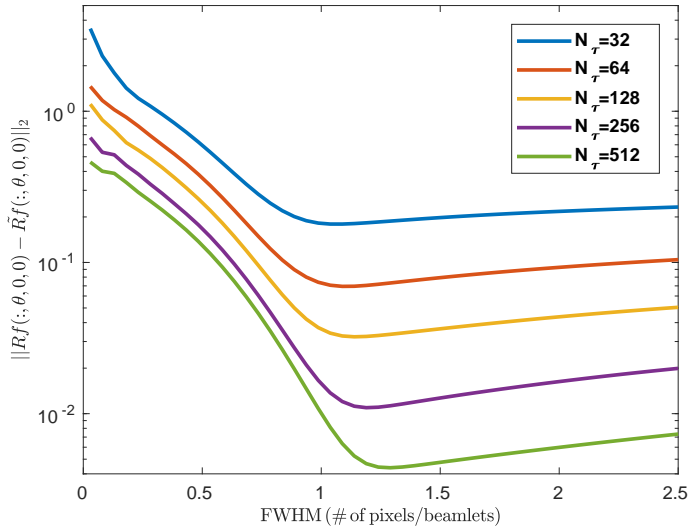}
  \caption{Approximation error incurred in the computation of the translation 
$Rf(\tau - P_\theta)$ when replacing the Dirac delta by the Gaussian 
\eqref{EQN:Gaussian} and carrying out the translation via the discrete Fourier 
transform.  The error is shown as a function of the Gaussian parameter $\sigma$ 
and of $N_\tau$, the number of beamlets of data available.}
  \label{fig:sigmaError}
\end{figure}

\mdone{R1C2}\changed{We choose to have the full width at half maximum (FWHM) 
of the Gaussian regularizer to cover one unit of beamlet width. 
Since the FWHM for a Gaussian with variance $\sigma^2$ is given by ${\mathrm  {FWHM}}\approx 2.355\sigma$, this means that we choose 
$\sigma = 1/2.355 \approx 0.42$. Observe that this choice is supported by \Fig{fig:sigmaError}, which displays the approximation error in Eqn.~\eqref{EQN:TranslationError}
in the simulated beam data from 
\Fig{fig:p_phantom}. 
For the values of $N_\tau$ considered, the choice ${\mathrm  {FWHM}}=1$ is
roughly the point at which the second error term in
\eqref{EQN:TranslationError} takes over from the first.}


\subsection{Numerical Illustration}
\label{SSEC:Numerical}

\begin{figure}[t]
  \centering
  \subfloat[]{\label{fig:p1_phantom}\includegraphics[width=0.6\textwidth]{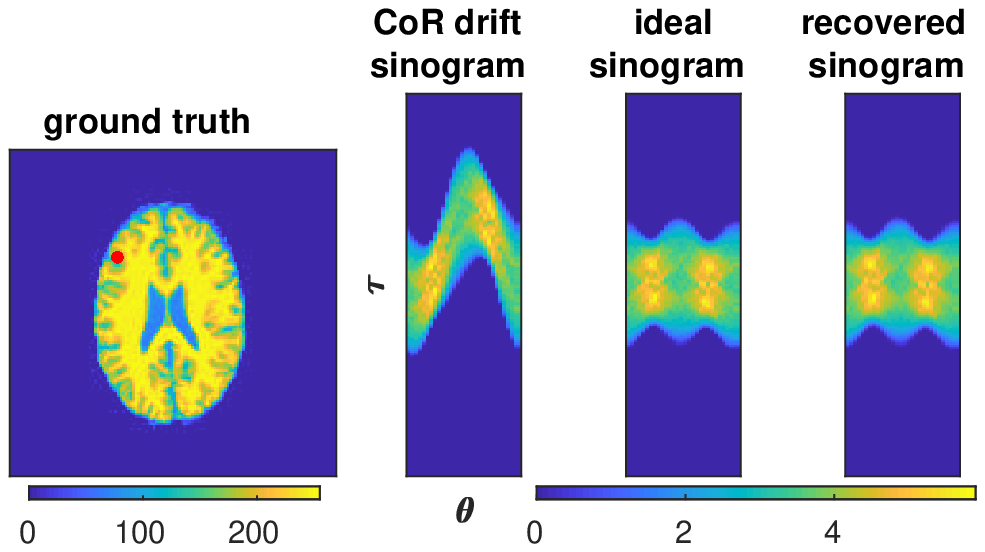}} \\
  \subfloat[]{\label{fig:pm_phantom}\includegraphics[width=0.6\textwidth]{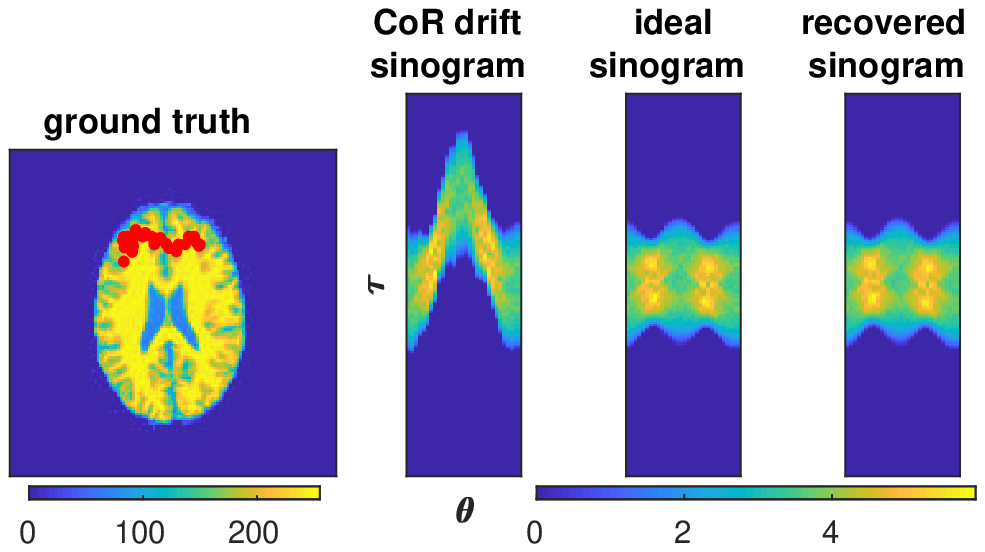}}
  \caption{Simulated numerical illustration of recovering the sinogram using the procedure of Sec, \ref{SSEC:Aligning}.  (a) Recovering a single CoR shift, made once at the beginning of the experiment.  (b) Recovering multiple CoR shifts, one at each angle of rotation.  The red dots in the ``ground truth'' images show the path of the CoR as it drifts during the data acquisition process.}
  \label{fig:p_phantom}
\end{figure}

To illustrate the effects of the approximation \eqref{EQN:TranslationError} in
the context of a complete image, we perform a simulation using a standard MRI
test image.  We simulate beam data from a variety of angles, assuming that
as the object rotates, it does so about a CoR 
$(x^*_{\theta_i},y^*_{\theta_i})= (x^*_{\theta_j},y^*_{\theta_j}) \neq (0,0)$ where $\theta_i \in 1,\dots,N_\theta$ and $\theta_j \in 1,\dots,N_\theta$ represent two different angles; in other words, there is only one CoR drift in the beginning of the experiment.  The
results are displayed in \Fig{fig:p1_phantom}, where $N=128$, $N_\tau=\lfloor \sqrt{2}N \rfloor$, and $N_\theta=30$ equally spaced in $[0,2\pi)$. 
The leftmost panel shows
the test image; the red dot 
marks the 
CoR $(x_\theta^*, y_\theta^*)$.  The second panel shows
the sinogram obtained from this data.  The third panel shows the sinogram that
we would have obtained in the absence of the drift in CoR.  The fourth
panel shows the sinogram that we obtain after performing the recovering
translations using the approach outlined in the preceding subsection.  The
drift-free and drift-recovered sinograms are in excellent agreement, even in
the presence of the errors described by \eqref{EQN:TranslationError}. 

In \Fig{fig:pm_phantom}, we repeat the experiment of
\Fig{fig:p1_phantom} except that the CoR changes each time the object is
rotated instead of just once at the beginning of the experiment; in other words, $(x^*_{\theta_i},y^*_{\theta_i})\neq (x^*_{\theta_j},y^*_{\theta_j}) \neq (0,0)$. The trail of
red dots in the leftmost panel shows how the CoR drifts across the image as the
data is acquired.  The drift-free and drift-recovered sinograms are again in
excellent agreement.  From a practical perspective, the additional ``jitter''
due to the multiple shifts of CoR---readily apparent in the measured
sinogram---complicates the alignment process significantly.  Nevertheless, the
algorithm we propose next is able to handle multiple
shifts.

\section{Optimization-Based Reconstruction Algorithm}
\label{sec:alg}

In the experiments of \Sec{SSEC:Numerical}, we assumed exact 
knowledge of the CoR
$(x_\theta^*, y_\theta^*)$ and the corresponding translation parameter $P_\theta$.  In 
practice, we must
estimate these parameters from the data.  As described in
\Sec{sec:intro}, the best techniques used to date generally do this
iteratively by alternating between using the current parameters $P_\theta$ to
reconstruct the image and using the reconstruction to update the 
parameters. Instead, we propose to estimate the parameters and
reconstruct the image concurrently.

The most straightforward approach would be to try to estimate $x_\theta^*$ and $y_\theta^*$
and obtain $P_\theta$ from these estimates using \eqref{EQN:RTFixCoR2}; we
refer to this as the \emph{explicit} approach, since we try to find the CoR
explicitly.  Alternatively, we can skip the estimation of $x_\theta^*$ and $y_\theta^*$ and
try to solve for $P_\theta$ directly; we call this the \emph{implicit}
approach. One major difference between the explicit and implicit approaches is the different dimension of unknown variables. In general, each $P_\theta$ corresponds to a pair of $(x_\theta^*,y_\theta^*)$, which results in $N_\theta$ number of CoRs to be recovered. If one knows which angles have the CoR drifts, however, the number of $(x_\theta^*,y_\theta^*)$ pairs can be reduced to the exact number of drifts.  

To exploit the correlations between $(x_\theta^*,y_\theta^*)$ 
and
$\bcW$, we formulate the reconstruction problem as simultaneously recovering 
the CoR and recovering the object. Therefore, the final \emph{explicit} reconstruction
problem is
\begin{equation} \label{1cor_obj} \ds
  \min_{\bcW\geq 0,\bstar{x},\bstar{y}}\phi(\bcW,\bstar{x},\bstar{y})= \frac{1}{2}\left|\left| \Le
  \bcW-\mathrm{vec}\left(\bg(\D,\bstar{x},\bstar{y})\right) \right|\right|^2_2,
\end{equation}
where $\bstar{x}=[x^*_\theta]_{\theta=1}^{N_\theta}$, $\bstar{y}=[y^*_\theta]_{\theta=1}^{N_\theta}$, $\bcW\geq 0$ is due to the physical nature of mass, $\D\in R^{N_\theta\times N_\tau}$ is the measurement data,
$\bg(\D,\bstar{x},\bstar{y})=\left[ D_{\theta,\tau}*\delta_{P_\theta(x_\theta^*,y_\theta^*),\sigma}\right ]_{\theta,\tau} \in \R^{N_\theta \times N_\tau}$ is the 
translated sinogram by \eqref{EQN:RTFixCoR2}, \ssout{$\bcW\in \R^{N^2}$ 
is the vectorized image based on a $N\times
N$ 2D discretization,} and $\Le=[L_v^{\tau,\theta}]\in \R^{N_\theta N_\tau\times N^2}$ is 
determined based on the standard 2D Radon mapping with an implied 
CoR of $(0,0)$.

As a side note, tomographic reconstruction without CoR error
recovery is typically formulated as
\begin{equation} \label{tomo_obj}
\ds \min_{\bcW\geq
0} \frac{1}{2}\left|\left| \Le \bcW-\mathrm{vec}\left(\D\right) \right|\right|^2_2,
\end{equation}
\added{which is equivalent to Eq.~\eqref{1cor_obj} for $(\bstar{x},\bstar{y}) = 0.$} 

The first-order
derivative of the objective function~\eqref{1cor_obj} is
\begin{equation}
  \label{D1cor} 
  \begin{array}{rcl}
   \nabla\phi (\bcW,\bstar{x},\bstar{y}) &=& \left[ \begin{array}{l} 
    \ds \nabla_{\bcW} \phi(\bcW,\bstar{x},\bstar{y}) \\ 
    \ds \nabla_{\bstar{x}} \phi(\bcW,\bstar{x},\bstar{y}) \\ 
    \ds \nabla_{\bstar{y}} \phi(\bcW,\bstar{x},\bstar{y}) 
  \end{array} \right]
\\  &=&
  \left[ \begin{array}{ll} \ds
    \Le^T \\ \ds
    \left(\nabla_{\bstar{x}}g(\D,\bstar{x},\bstar{y})\right)^T \\ \ds 
    \left(\nabla_{\bstar{y}} g(\D,\bstar{x},\bstar{y})\right)^T \\
  \end{array} \right] \left(\Le \bcW-\bg(\D,\bstar{x},\bstar{y}) \right).
  \end{array}
\end{equation}
\added{In \eqref{D1cor}, $\nabla_{\bstar{x}}g(\D,\bstar{x},\bstar{y})$ and
$\nabla_{\bstar{y}} g(\D,\bstar{x},\bstar{y})$ are}
\[
\ds \mathrm{vec}\left(\left[\D_{\theta,:}*\left(\frac{1}{\sigma \sqrt{2 \pi}}\exp \left(\frac{-(\cT-P_\theta)^2}{2\sigma^2}\right)\right)\circ \frac{\cT-P_\theta}{\sigma^2}\left(\cos\theta-1\right)\right]_{\theta=1,\dots,N_\theta}\right),
\]
and
\[
\ds \mathrm{vec}\left(\left[\D_{\theta,:}*\left(\frac{1}{\sigma \sqrt{2 \pi}}\exp \left(\frac{-({\cT}-P_\theta)^2}{2\sigma^2}\right)\right)\circ \frac{{\cT}-P_\theta}{\sigma^2}\sin\theta\right]_{\theta=1,\dots,N_\theta}\right),
\]
\added{respectively,} where $\circ$ is the Hadamard product, and $\D_{\theta,:}$ is the $\theta\,th$ row of $\D$.
\ssout{To recover for the multiple CoR-shifts case, we extend the objective
function in \eqref{1cor_obj} to multiple $(x^*_\theta,y^*_\theta)$
for each $\theta$.} \mdone{R1C3} A major difference compared with the single CoR-shift case 
is the 
increased computational cost. Instead of having only two extra
parameters, this case has $2N_\theta$ extra parameters.

For now, we have been focusing on explicitly recovering the 
coordinates of the CoRs $(x_\theta^*,y_\theta^*)$. Alternatively, we can reformulate the
optimization problem as finding the optimal shifts $P_\theta$ for each angle, 
which results in the \emph{implicit} problem 
\begin{equation}\label{eqn:implicit}
\ds \min_{\bcW\geq
  0,\cP} \phi(\bcW,\cP)=\frac{1}{2}\left|\left| \Le \bcW-\bg(\D,\cP)
  \right|\right|^2_2,
\end{equation}
where $\cP=[P_\theta]_{\theta=1,\dots,N_\theta}$ and $\bg(\D,\cP)=\left[ D_{\theta,\tau}*\delta_{P_\theta,\sigma}\right ]_{\theta,\tau} \in \R^{N_\theta \times N_\tau}$. 
The derivative of the 
objective function in \eqref{eqn:implicit} is similar to 
\eqref{D1cor} but simpler, so we
will not describe it here. In the case of multiple CoR shifts, the advantage of
this formulation is the reduced number of extra parameters by a factor
of 2.
Notice that given a fixed CoR, the Radon transform has the property that different transformations (e.g., translation or rotation) of an object will result in different sinograms;  see Appendix \ref{sec:properties_radon}. These properties suggest that any pair of transformed object and its corresponding sinogram is an optimal solution of the reconstruction. 
Therefore, the
proposed optimization problem, either explicit or implicit, 
will not have a
unique solution. In other words, if two objects are related by equal (up to 
translations and rotations), then the two objects, together with their 
corresponding
sinograms, are equivalent from the perspective of the 
optimization problems in \eqref{1cor_obj}
and \eqref{eqn:implicit}.


\section{Optimization Complexity and Computational Expense}
\label{sec:opt}

\begin{figure}[!htbp]
\centering
\includegraphics[scale=0.7]{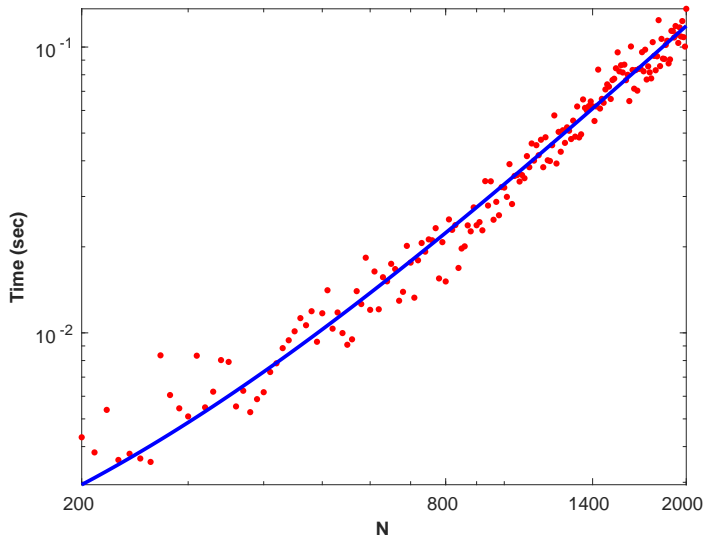}
\caption{Time elapsed for one (function, gradient) evaluation for different 
\changed{measurement resolution numbers} 
 $N$.}
\label{fig:time}
\end{figure}

Since the complexity of Eq.~\eqref{EQN:RTFixCoR2} is negligible, we analyze only the computational complexity of the implicit approach~\eqref{eqn:implicit}. The calculation of Eq.~\eqref{eqn:implicit} requires about $N^2N_\theta$ flops given the relationship $N_\tau=\lfloor\sqrt{2}N\rfloor$, and it includes 1 misfit calculation and 3 FFTs required in operator $\bg$. \Figure{fig:time} shows \added{a log-log plot of} \mdone{R1C4} the computational time of one (function, gradient) evaluation for increased number of $N$, given $N_\theta=1$. With a model fit to the time result, the time complexity is on the order of $N^2$ which is consistent with our analytical approximation considering $N_\theta \ll N$ in general.

The optimization solver we use in the numerical experiments is an
inexact truncated-Newton (TN) method \cite{Nash200045} with a preconditioned
conjugate gradient (PPCG) 
\cite[pp.~460]{nocedal2006numerical} subproblem solver to compute the search 
direction. The  result, summarized in Alg.~\ref{TN}, is a large-scale 
optimization solver with memory
efficiency well suited for the high computational complexity required by the
resulting nonlinear, nonconvex optimization problem.

\begin{algorithm} \caption{Truncated Newton algorithm for the implicit 
problem \eqref{eqn:implicit}.} 
   \label{TN}
\begin{algorithmic}[1]
\STATE Input $\left (\bcW^{(0)},\cP^{(0)} \right)$ and tolerance $\epsilon>0$; 
set $k=0$.
 \REPEAT
 \STATE Obtain search direction $d^{(k)} \leftarrow 
\text{PPCG}\left(\left(\bcW^{(k)},\cP^{(k)}\right),\nabla^2 \phi,\nabla 
\phi,\mathbf{0}\right).$ If no descent direction is obtained, switch to steepest descent.
\STATE Obtain $\alpha^{(k)} \leftarrow
\text{Projected Line 
Search}\left(d^{(k)},\left(\bcW^{(k)},\cP^{(k)}\right),\phi,\nabla 
\phi,\mathbf{0}\right)$; see
\cite{lin1999newton}.\

\STATE Update 
$\left(\bcW^{(k+1)},\cP^{(k+1)}\right) \leftarrow
 \left(\bcW^{(k)},\cP^{(k)}\right)+\alpha^ { (k) } d^{(k)}$.
\UNTIL  The stopping criterion $\left| \left| \nabla 
\phi\left( \bcW^{(k)},\cP^{(k)} \right) \right|
\right|\leq \epsilon$ is satisfied.
\end{algorithmic}
\end{algorithm}

In our implementation, we do not form the Hessian $\nabla^2 \phi$; 
instead, we approximate the Hessian-vector product $(\nabla^2 \phi)^T d$ 
required in PPCG by taking finite differences with $\nabla \phi$ values. 
An estimation of the complexity of TN is provided as follows. Each outer
iteration needs the following computations with respect to the number of
unknown parameters:
\begin{itemize}
  \item 1 infinity-norm calculation, 1 vector addition, and 2
    (function, gradient) evaluations
  \item a number of PPCG iterations, with cost per inner iteration given by
  \begin{itemize}
    \item 1 (function, gradient) evaluation,
    \item 4 inner products, and
    \item 5 ``vector+constant$\cdot$vector'' operations
  \end{itemize}
\end{itemize}

In our experience, on average, 5 PPCG iterations are required per outer TN
iteration. The overall cost of one TN iteration for solving 
\eqref{eqn:implicit} amounts to $(1+5+20+25+1)(|\cV|+N_\theta)=52(|\cV|+N_\theta)$ floating-point operations, plus 7 (function, gradient) evaluations, whose complexity is $O(N^2)$.

\section{Numerical Results}
\label{sec:numerical}

In this section, we examine the performance of the algorithm for both the
explicit and implicit cases. The primary goal of our tests is to measure the
performance of how well we recover the CoR shift with respect to different initializations,
different objects, and different levels of noise in the data. All numerical
experiments are performed on a platform with 32 GB of RAM and two Intel
E5430 Xeon CPUs. Throughout all the numerical tests, we fix the experimental
setup for all the tests as $N=128$, 
$N_\tau=\lfloor\sqrt{2}N\rfloor$, and
$N_\theta=30$; \mdone{R2C1}\added{that is, the object size is $128\times 128$ 
pixels, and the tomographic data is 
\changed{measured} 
by collecting 30 projections over a full rotation of $2\pi$ radians}. 
We also choose the stopping criteria of TN to be $\left|
\left| \nabla\phi \right| \right|<10^{-5}$.

As indicated in \Sec{sec:alg}, our proposed optimization
problem does not have a unique solution in the sense that the reconstructed
object can be a translated or rotated version of the ground truth. Therefore,
the error metric we use to measure the reconstruction quality cannot be
simply the mean squared error. To resemble the human 
visual system, we utilize the structural
similarity (SSIM) metric from \cite{wang2004image} to quantify the
reconstruction quality. Given two
images $a$ and $b$ of the same dimension, the SSIM index is a 
measure of  the similarity
between $a$ and $b$ and is defined by \[\ds \text{SSIM}(a,b)=\frac 
{(2\mu_{a}\mu _{b}+c_1) 
(2\sigma_{ab}+c_2)}{(\mu_{a}^{2}+\mu_{b}^{2}+c_1)(\sigma 
_{a}^{2}+\sigma_{b}^{2} + c_{2})}, \] where $c_1$ and
$c_2$ are small, positive constants; $\mu_a$ and $\sigma_a^2$ are the mean  
and variance, respectively, of $a$; $\mu_b$ and $\sigma_b^2$ are the mean  
and variance, respectively, of $b$; and 
 $\sigma_{ab}$ is the covariance of $a$ and $b$. Notice that for nonnegative 
$a$ and $b$, $\text{SSIM}(a,b) \geq 0$. The closer
the value of $\text{SSIM}(a,b)$ is to $1$, the greater is the similarity between $a$ and $b$.

\begin{figure}[!htbp]
\centering
\includegraphics[scale=0.8]{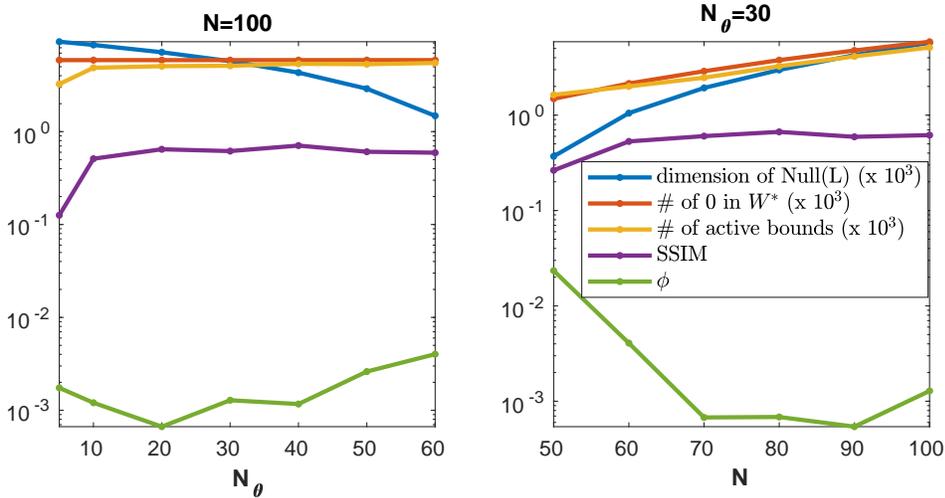}
\caption{Performance of TN for solving problem \eqref{tomo_obj} for various underdetermined and overdetermined systems.}
\label{fig:null}
\end{figure}

We emphasize that in this work our focus is on the development of a joint
framework to account for general CoR shifts. To isolate the contributions of our 
proposed approach to the reconstruction performance, we do not consider any 
\changed{additional} 
regularization technique in this work. Additionally, the 
nonnegative constraint on $\bcW$ serves as a soft regularizer, which means we do 
not need to add a regularizer and can instead show the benefit of our approach 
by itself \added{\cite{donoho2005sparse,bruckstein2008uniqueness}}. As suggested by \Fig{fig:null}, where we test the performance of TN 
on various settings for \added{the} traditional tomography problem 
\eqref{tomo_obj}, we observe that as long as the dimension of Null($\Le$) is 
close to the number of zeros \changed{of the object}, the reconstruction quality is 
stable and satisfactory based on SSIM and the objective value of problem 
(\ref{tomo_obj}). Since the object domain always has a zero support to account 
for the field of view during the rotation, the dimension of Null($\Le$) 
typically is proportional to the number of zeros in $\bcW^*$, and this explains 
the satisfactory reconstruction result from seemingly underdetermined system.  

\subsection{Analysis of Initialization}
\label{subsec:init}

In this section, we test the sensitivity of the 
proposed algorithm to the initialization of $\bcW$. 

\begin{figure}[!htbp]
\centering
\includegraphics[scale=0.8]{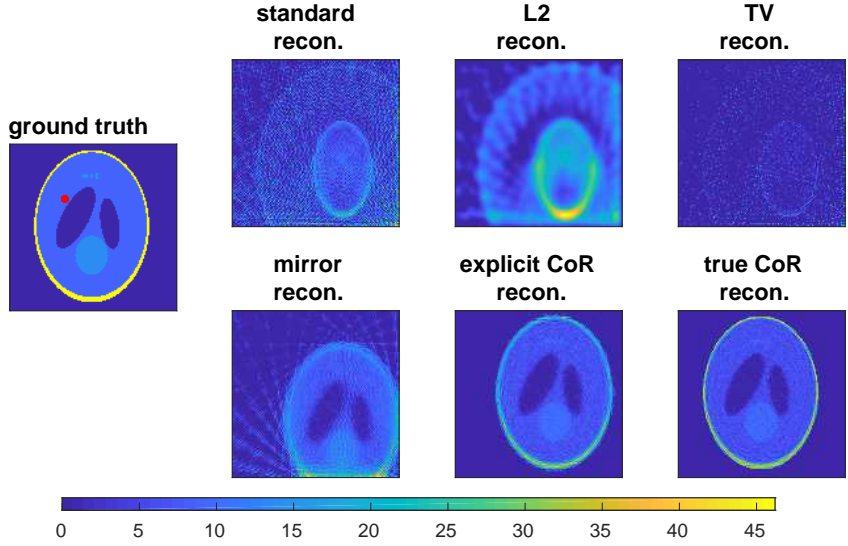}
\caption{(First) Shepp-Logan phantom object with true CoR  denoted by the red 
dot. \changed{(Upper left) Without error recovery, the direct reconstruction result is far 
away from the ground truth. (Upper middle) Result of  direct reconstruction with L2-regularization. \changed{The regularization parameters $\lambda$ for L2 and TV regularizers are chosen to be the maximizers of SSIM}. (Upper right) Result of direct reconstruction with TV-regularization.  (Lower left) Result of reconstruction by ``mirror alignment.'' (Lower middle)} Reconstruction with \added{explicit} CoR recovery. (Lower right) Result of reconstruction using the true CoR.}
\label{fig:1cor_x}
\end{figure}

\begin{figure}[!htbp]
\centering
\includegraphics[scale=0.7]{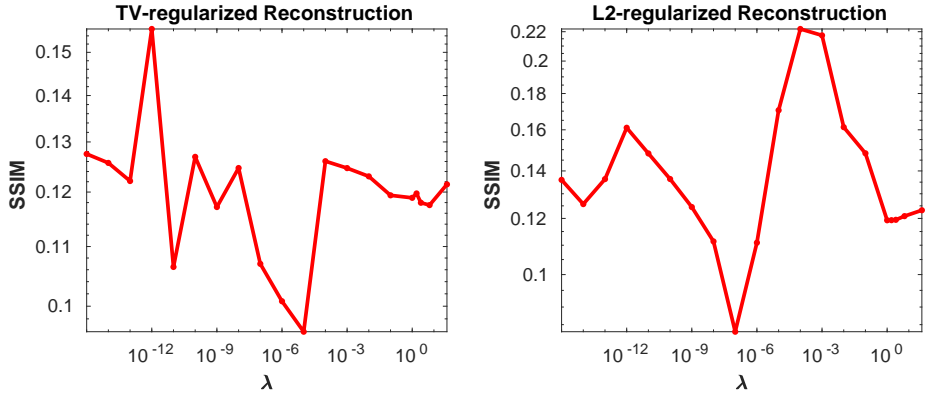}
\caption{ (Left) SSIM for TV-regularized reconstruction with different $\lambda$. (Right) SSIM for L2-regularized reconstruction with different $\lambda$.}
\label{fig:lambda}
\end{figure}

We first consider the case of a single 
CoR shift that happens at the beginning
of the experiment on a standard Shepp-Logan phantom object. The CoR 
is denoted by the red dot on top of the ground truth in
\Fig{fig:1cor_x} \added{(far left)}. \changed{The upper left panel} of \Fig{fig:1cor_x} 
shows the standard reconstruction without considering the shift of the CoR, and 
we observe that it is far from the ground truth. We also include reconstructions 
obtained by adding standard regularizers.
\changed{In particular, we consider $L2$- and $TV$-regularizers
\cite{li2016s, hansen2011total} on the variable $\bcW$, respectively given by}
\begin{equation} \label{eqn:L2}
\ds \min_{\bcW\geq
0} \frac{1}{2}\left|\left| \Le \bcW-\mathrm{vec}\left(\D\right) 
\right|\right|^2_2 + \lambda \left| \left| \mathcal{L}\bcW  \right| \right|_2^2
\end{equation}
\added{and}
\begin{equation} \label{eqn:TV}
\ds \min_{\bcW\geq
0} \frac{1}{2}\left|\left| \Le \bcW-\mathrm{vec}\left(\D\right) \right|\right|^2_2 + \lambda \left| \left| \nabla\bcW  \right| \right|_1,
\end{equation}
where $\mathcal{L}\in \R^{N^2\times N^2}$ is the Laplacian operator, \added{$\nabla$ denotes the spatial gradient}, and $\lambda>0$ is the regularizer parameter that balances the misfit term and regularization term. One approach to choose  $\lambda$ would be the L-curve method \cite{hansen1993use}. Instead, however, we choose $\lambda \in [10^{-10},10^2]$ to be the value that gives the best reconstruction as measured by the SSIM index (see \Fig{fig:lambda}). This is an idealistic choice of $\lambda$ that is impractical. The results given in the \changed{upper middle and upper right panels} demonstrate this assertion and show that standard regularization techniques are not able to improve the reconstruction quality in the presence of even a single CoR shift. Since even an optimal $\lambda$ performs poorly for CoR recovery, for the rest of our experiments we do not report results from regularization techniques. 

\mdone{R1C5}\added{We also compare our reconstruction result with the most popular approach for the case of a single CoR drift \cite{min2012new}.  The basic idea is to utilize the fact that when the CoR is the origin, the projection at angle $0^\circ$ should be a reflection of the projection at angle $180^\circ$. Therefore, the CoR shift can be estimated by aligning these two projections. In our case, as is often the case in practice, exactly reflecting directions are not available; instead, we choose the pair of angles $1^\circ$ and $175^\circ$ from the available 30 projections. We align these two projections by cross-correlation \cite{zitova2003image} and use the estimated CoR shift to translate the collected sinogram. The reconstruction from the shifted sinogram is reported in the lower left panel. We also report the reconstruction using the true CoR in the lower right panel. The lower middle panel shows the reconstruction obtained by the proposed explicit approach, which clearly resembles the ground truth the best.}

\begin{figure}[!htbp]
\centering
\includegraphics[scale=0.8]{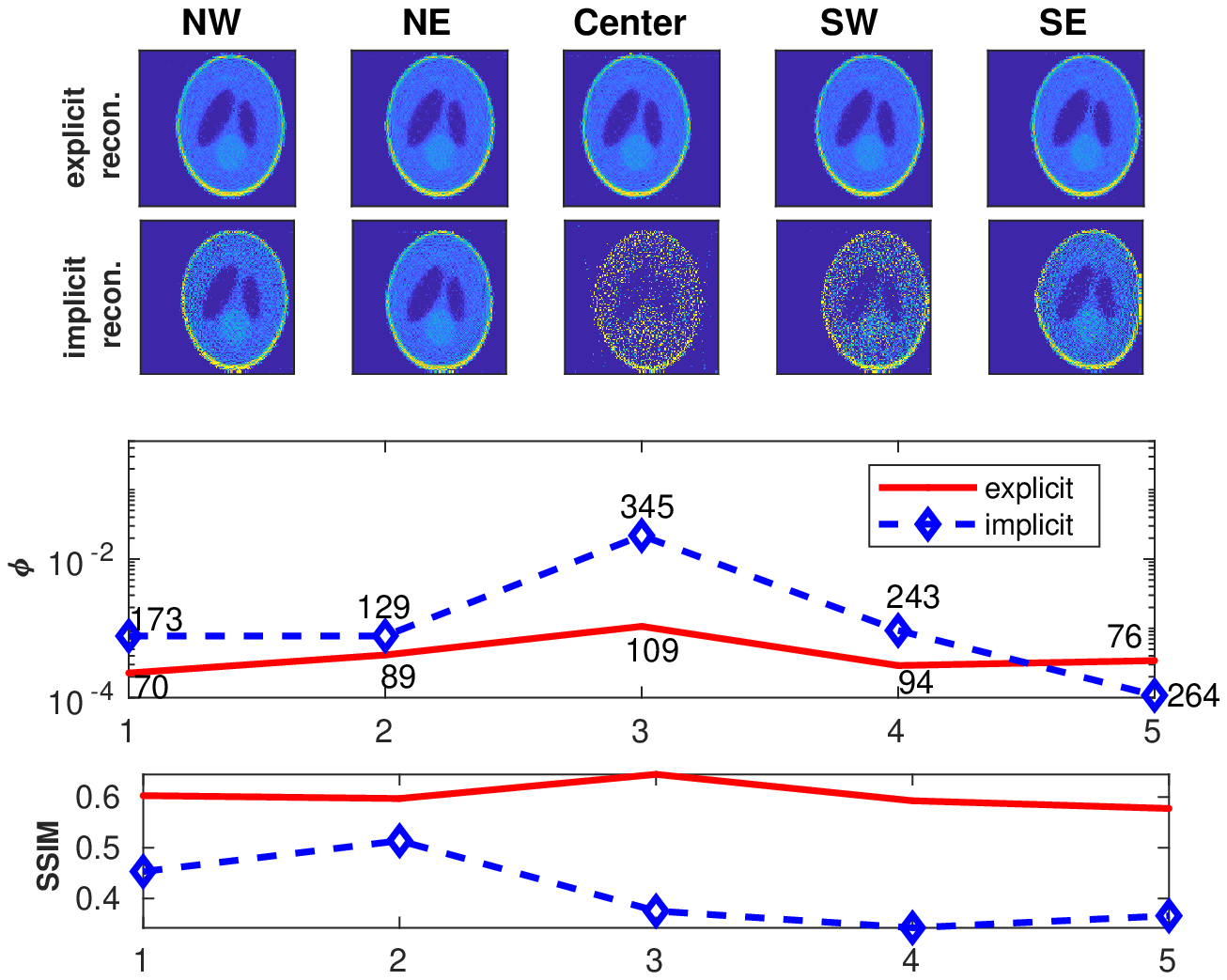}%
\caption{Performance of the explicit and implicit approaches for recovering single 
CoR shift. Columns correspond to different initializations for CoR 
$(x_\theta^*,y_\theta^*)$. The top row shows the corresponding reconstructions from the explicit approach; the second row shows the corresponding 
reconstructions from the implicit approach; the third row shows the 
corresponding objective function values from the explicit and implicit 
approaches, respectively, where the number labeled on each case is the total 
number of iterations needed for convergence; the fourth row shows the 
reconstruction qualities measured by SSIM index for the explicit and implicit 
approaches, respectively.}
\label{fig:1cor_ini}
\end{figure}

The \changed{next} test is to optimize the
explicit problem~\eqref{1cor_obj} and the implicit problem~\eqref{eqn:implicit}
in order to recover $(x_\theta^*,y_\theta^*)$ and $\cP$, respectively. Given an initial
 $(x_\theta^*,y_\theta^*)$ as the center, left top corner (NW), left
bottom corner (SW), right top corner (NE), and right bottom corner (SE) of the
object domain, respectively, and initialization $\bcW^*=0$, the performance
is shown in \Fig{fig:1cor_ini}.  In this test, the explicit approach 
consistently recovers objects that are better compared with the implicit approach, which struggles 
for certain initializations. This observation is also reflected in the 
objective value and the
reconstruction quality. In addition, the explicit approach takes fewer iterations
to converge to smaller objective values compared with the implicit approach.
The explicit approach also provides relatively better
reconstruction quality as measured by SSIM index.

\begin{figure}[!htbp]
\centering
\includegraphics[scale=0.8]{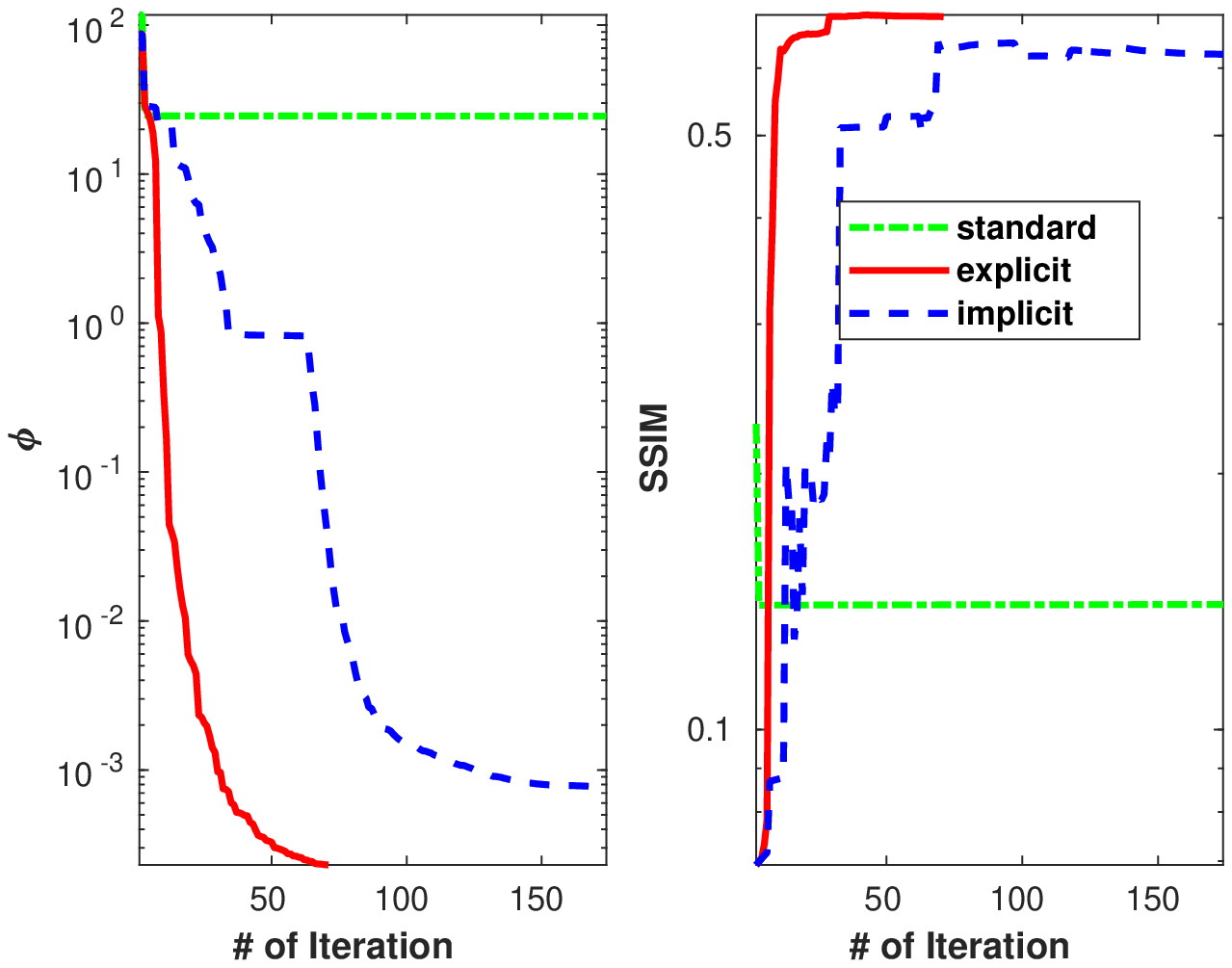}%
\caption{Iterative performance of three approaches---standard tomographic reconstruction, explicit reconstruction, and implicit reconstruction---for the single CoR-shift case. Left: progress of reducing objective value $\phi$. Right: progress of improving SSIM.}
\label{fig:1cor_iter}
\end{figure}

In \Fig{fig:1cor_iter}, we also show the iterative progress of the
optimization problems \eqref{tomo_obj}, \eqref{1cor_obj}, and \eqref{eqn:implicit}
respectively for the single CoR-shift case. The left graph shows the progress of the objective value for the three approaches. Standard tomographic
reconstruction, which refers to reconstruction without error recovery,
converges much more slowly than the other two with a much higher objective value,
while the explicit approach reduces the objective value even further than
the implicit approach because of the reduced number of variables. The right graph compares the progress of the SSIM index for the three approaches, which is somewhat
consistent with the function value reduction. We observe that the explicit approach
provides the best reconstruction quality in terms of SSIM compared with the other two
approaches, while the standard approach performs the worst because it does not attempt to recover for the error in CoR.

\begin{figure}[!htbp]
\centering
\includegraphics[scale=0.8]{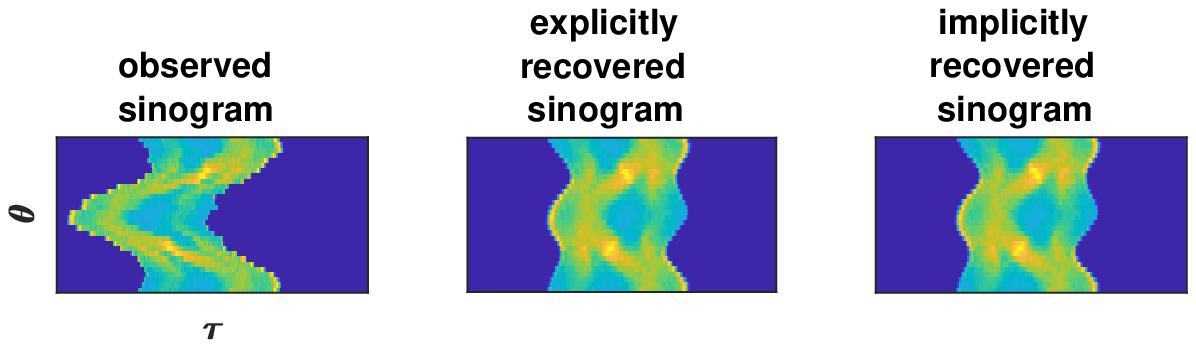}%
\caption{(Left) Multiple CoR shifts introduce jitters in the measurement data. (Middle) the explicit approach provides a more clean and jitter-free sinogram. (Right) The implicit approach is also able to recover a clean and jitter-free sinogram. }
\label{fig:sino}
\end{figure}

\begin{figure}[!htbp]
\centering
\includegraphics[scale=0.8]{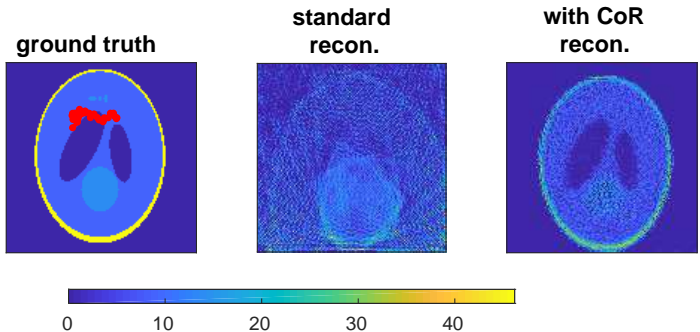}
\caption{(Left) Shepp-Logan phantom 
\changed{object} 
with true CoRs connected by red dots. (Middle) Without 
error recovery, a direct reconstruction result is far away from the ground 
truth. (Right) Reconstruction result with error recovery, which is very close to 
ground truth.}
\label{fig:mcor_x}
\end{figure}

\begin{figure}[!htbp]
\centering
\includegraphics[scale=0.8]{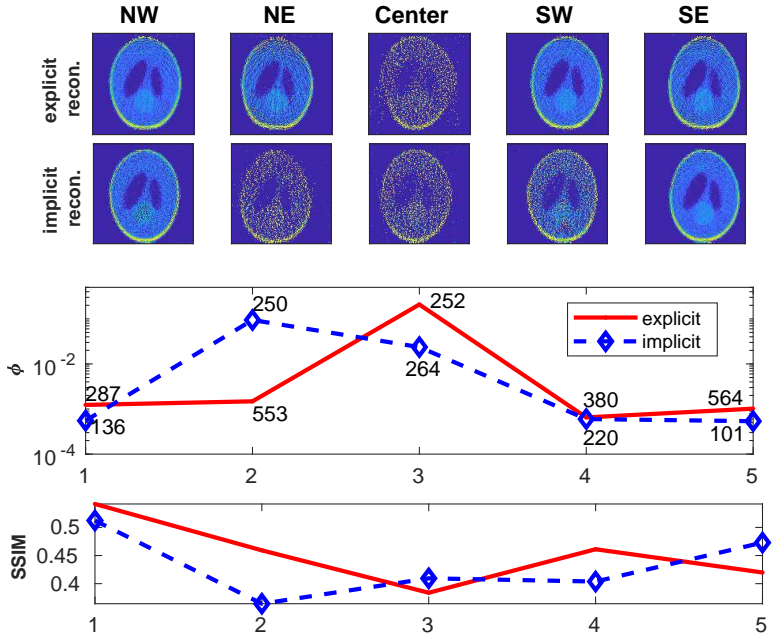}%
\caption{For the multiple CoR-shifts case, the reconstruction qualities provided by the explicit and the implicit approaches are comparable, but the implicit approach converges faster. }
\label{fig:mcor_ini}
\end{figure}

The next test concerns solving the reconstruction problem with multiple CoR
shifts in the sense that every rotation has its own CoR. This problem setting
is illustrated in \Fig{fig:mcor_x}, where the ground truth is again the
Shepp-Logan phantom and the trajectory of the dynamic CoR is 
denoted by connected red dots. The second panel shows the reconstruction without 
CoR
recovery, and the third panel shows the reconstruction with error recovery. Again, we tested the performance of the explicit and the
implicit algorithm; the results are shown in
\Fig{fig:mcor_ini}.  This time we observe that the reconstruction quality of the
implicit approach is comparable to the explicit approach but is more sensitive to the initialization. 
In the multiple CoR-shifts case, however, the convergence of the implicit approach
is faster because of its reduced number of variables compared with explicit approach.
Again, the performance difference between the explicit and implicit approaches
is reflected in the objective value and the corresponding SSIM index. More importantly, as shown in \Fig{fig:sino},
both approaches are able to remove the jitters introduced by multiple CoR
shifts and return a clean sinogram. The consistency between the objective
value $\phi$ and SSIM index suggests that, in the case of limited knowledge about the
object, the objective value can be reliable to judge the reconstruction quality of different local minima.

We note that as long as the shift parameters $\cP$ are initialized to zero, the performance is relatively stable and less susceptible to
the initialization of the object. The reason is that choosing the initial shift 
as
zeros is equivalent to providing the initialization of the object as the
reconstruction without error recovery, which is closer to the true solution than
the reconstruction from a randomly perturbed sinogram.  Therefore, choosing zeros
for the initial shifting parameters performs better, and this initialization will be our choice for the rest of the
tests.

\begin{figure}[!htbp]
\centering
\includegraphics[scale=0.8]{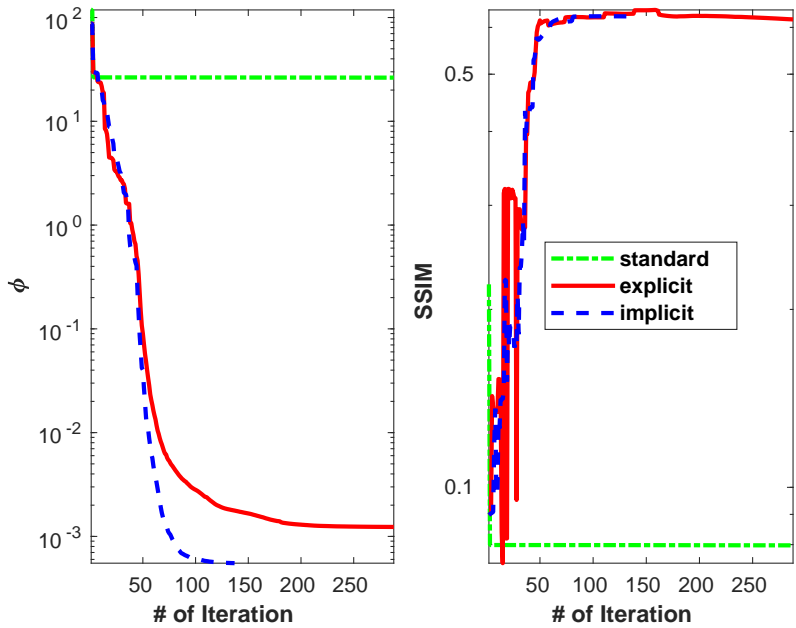}%
\caption{Iterative performance of three approaches---standard tomographic reconstruction, explicit reconstruction, and implicit reconstruction---for multiple CoR-shift case. Left: progress of reducing objective value $\phi$. Right: progress of improving SSIM.}
\label{fig:mcor_iter}
\end{figure}

Similar to the single CoR case, we compare the iterative progress provided
by the standard reconstruction, the explicit approach, and the implicit approach for the
multiple CoR-shifts case; the results are shown in \Fig{fig:mcor_iter}.
As
expected, the standard approach converges rapidly to a local minimum that 
has a high objective value. More importantly, the implicit approach provides a better objective value than does the explicit approach. The
reconstruction quality (right panel of \Fig{fig:mcor_iter}) shows 
comparable
performance between the explicit and the implicit approaches.

\subsection{Noisy Data}

\begin{figure}[!htbp]
\centering
\includegraphics[scale=0.8]{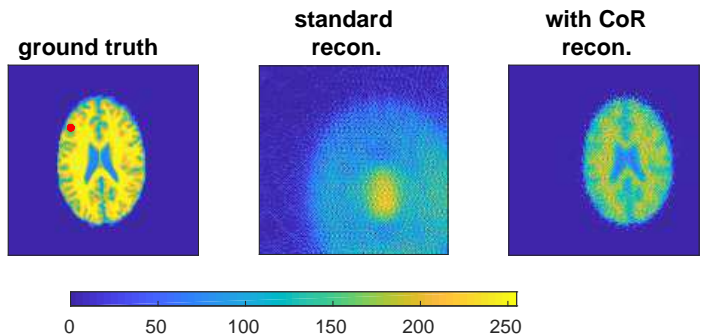}
\caption{(Left) Object as the standard 
MRI image with only one CoR shift denoted by the red dot. (Middle) 
Reconstruction from noise-free data without error recovery. (Right) Reconstruction with error recovery.}
\label{fig:1cor_phantom}
\end{figure}

In practice, the experimental data is known to be corrupted by
noise. In this section, we test the performance of our
proposed approach for contaminated data with respect to different levels of \added{(mean zero)  Gaussian}
noise. We start the test for the single CoR shift case as shown in
\Fig{fig:1cor_phantom}, where the ground truth in this case is the
standard MRI image. Similarly, the second and the third panels show the reconstructions without and with error recovery, respectively. 

\begin{figure}[!htbp]
\centering
\includegraphics[scale=0.8]{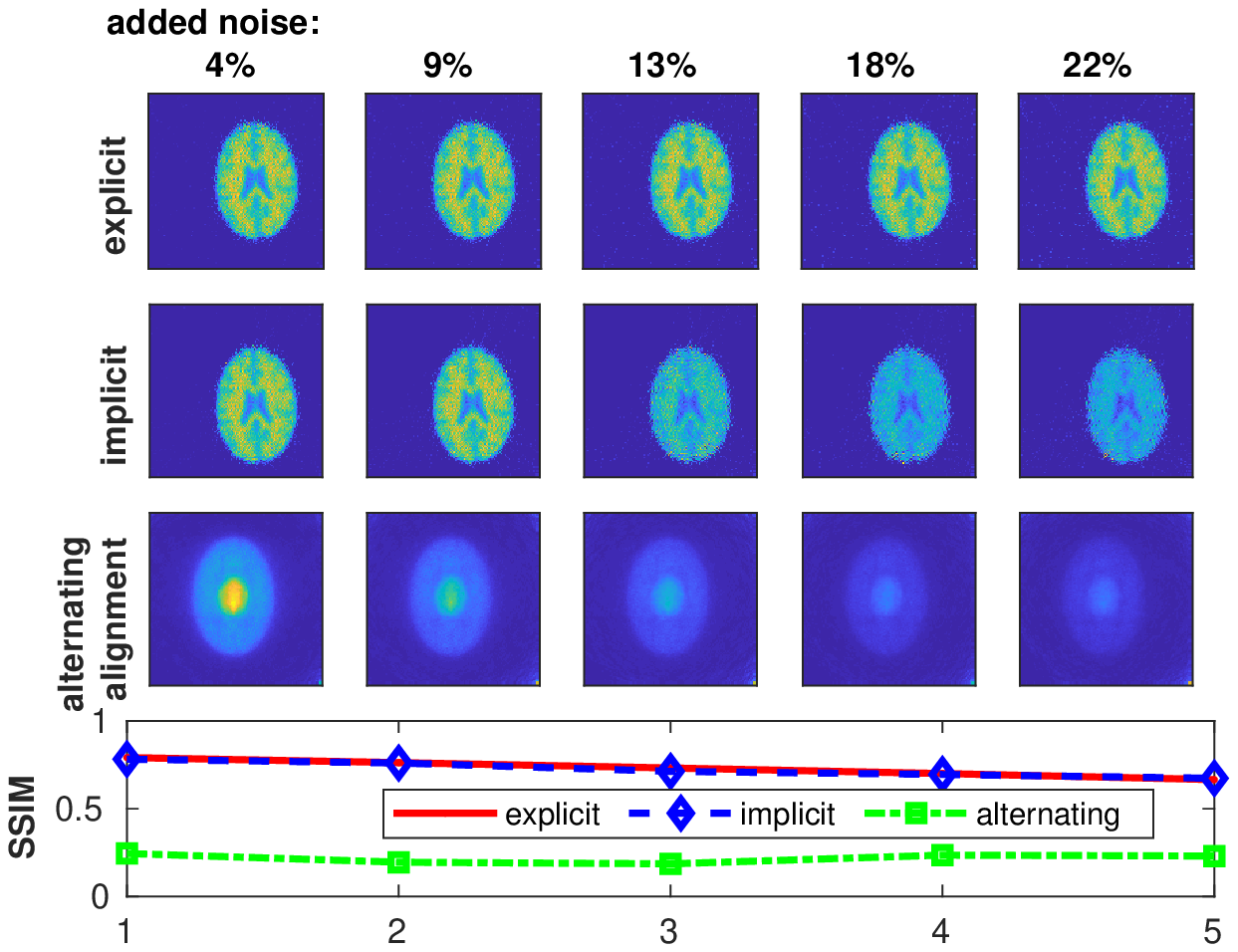}%
\caption{Performance comparison between proposed approach and state-of-the-art 
approach provided by TomoPy library for single CoR shifts when the noise level 
of experimental data is increasing.}
\label{fig:1cor_noise}
\end{figure}

In this case, we compare the performance of the explicit and implicit
approaches with the state-of-the-art alignment approach
\cite{gursoy2017rapid} provided by TomoPy 1.1.3 \cite{gursoy2014tomopy}, a
widely used tomographic data processing and image reconstruction library. 
The algorithm proposed in \cite{gursoy2017rapid} is an iterative
method that alternates between the alignment of the sinogram and the reconstruction
of the object. The alignment step minimizes the misfit between the
current sinogram and the newly simulated sinogram based on the current
reconstruction. In \Fig{fig:1cor_noise}, we illustrate the
reconstruction results returned by the three methods as we
increase the noise level in the experimental data; the columns
represent an increase in the noise level from 4\%
to 22\%. We observe that the explicit and implicit approaches are comparable
with each other. However, the SSIM index for both approaches suggests that the implicit approach provides better reconstruction quality. This observation can be explained by the fact that although the explicit approach can provide a higher contrast image, its ability to remove noise outside of the sample region is not as good as that of the implicit approach. Overall, the performance of our
proposed approaches outperforms the alternating algorithm
\cite{gursoy2017rapid} provided by TomoPy.

\begin{figure}[!htbp]
\centering
\includegraphics[scale=0.8]{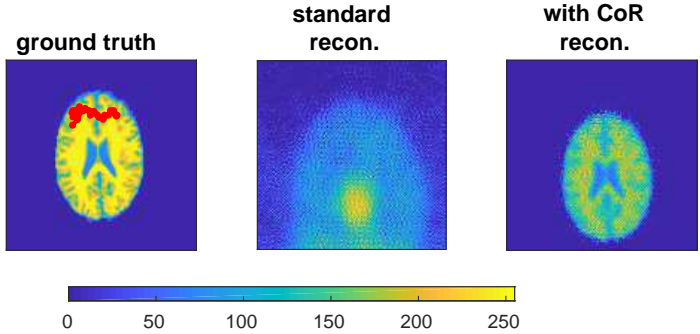}
\caption{(Left) Second test case given the MRI imaging object, with multiple CoR shifts with its trajectory connected by red dots. (Middle) Reconstruction from noise-free data without error recovery. (Right) Reconstruction of noisy data with error recovery.}
\label{fig:mcor_phantom}
\end{figure}

\begin{figure}[!htbp]
\centering
\includegraphics[scale=0.8]{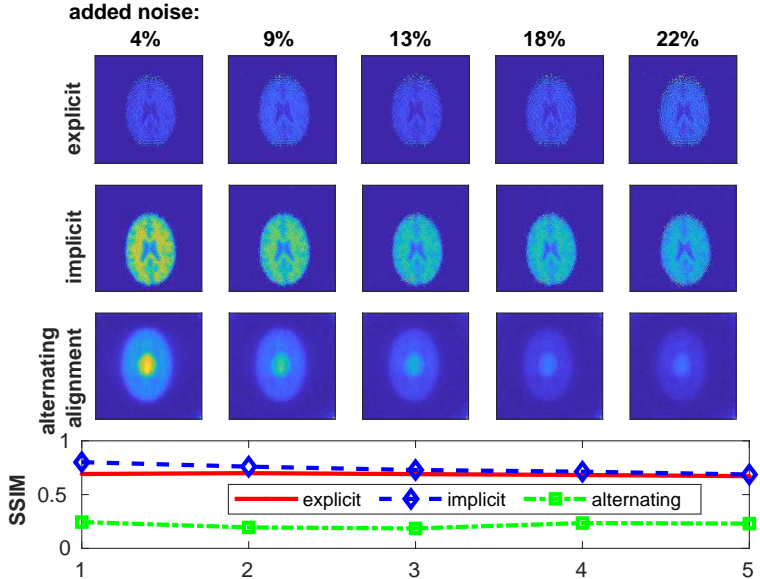}%
\caption{Performance comparison between proposed approach and state-of-the-art approach provided by TomoPy library for multiple CoR shifts when the noise level of experimental data is increasing.}
\label{fig:mcor_noise}
\end{figure}

The last study is to test the performance of our proposed approach for
recovering the multiple CoR shifts, as illustrated in
\Fig{fig:mcor_phantom}. As measured by SSIM, the implicit approach clearly outperforms the other two approaches while the reconstructions of the alternating alignment do not approach the ground truth at all. In \Fig{fig:mcor_noise}, we illustrate the
reconstruction results returned by the three methods respectively as we
increase the noise level in the same way as the single CoR shift case 
illustrated
earlier. Again, the explicit and implicit approaches are able to provide 
a reasonable reconstruction, and the implicit approach in this case is able to 
provide a better reconstruction for the noisier case. Similarly to the 
single CoR shift
case, the performance of our proposed approaches outperforms that of the alternating
algorithm dramatically. Therefore, in the case of limited knowledge of the CoR shifts with noisy measurement data, the implicit approach works the best. 

\section{Conclusions}
\label{sec:conc}
In this work, we propose a simultaneous inversion framework to
address a common yet challenging experimental error in tomographic
reconstruction, namely, a shift of the center of rotation. We also derive an analytical model to
describe the general center-of-rotation shifts and how the sinograms with and without such 
shifts relate to each other. We derive an analytical model for 
the reconstruction framework to recover the center-of-rotation shifts together
with the reconstruction of the object. The resulting optimization framework can
be tuned to explicitly or implicitly reveal the true CoR locations. It is
solved by a derivative-based optimization solver such as a  truncated Newton
algorithm.

Our numerical results show that when limited and
noisy experimental data is available, significant improvements of the
reconstruction quality can be achieved by performing the proposed simultaneous
inversion. We compare our new approach with a state-of-the-art error reconstruction technique implemented in TomoPy and standard reconstruction techniques with a regularizer approach. For both comparisons, we observe superior results obtained by our approach. Given prior knowledge of how center-of-rotation shifts happen, explicit and
implicit approaches can be chosen to maximize the performance, based on their
different computational complexities. In particular, the implicit approach
shows a stronger capability in more general scenarios such as noisy data. 

{\bf Discussion:}
The simultaneous inversion we propose in this work is
flexible and does not require extensive knowledge about the experimental system
or the object. \mdone{R1C1}\added{Even for the single CoR shift, the proposed approach does not require that the sample be completely contained in the field of view or that no nonlinear effects such as diffraction or refraction occur.} It can be naturally extended to partially address more complicated
experimental error such as per-scan drift. \mdone{R2C2}\added{In parallel-beam tomography, the trajectories of the penetrating x-rays result in tomographic projections that form a stack of parallel 2D slices. Traditionally, 3D objects are recovered by reconstructing each slice independently \cite{webb1990watching}. Therefore, the extension of the proposed approach to a full 3D reconstruction is trivial in parallel-beam tomography. Alternatively, since each 2D slice shares the same CoR, the proposed approach can also be extended to a direct 3D reconstruction without adding an extra degree of freedom.} One aspect we
can further explore is to incorporate prior knowledge about the
experimental system, such as the continuous dynamic of the center-of-rotation movement, into
the joint framework as a regularizer on the error parameters.  Overall, the
challenge we are facing is the increased computational complexity. One way to
mitigate this difficulty is to consider a parallel implementation that could
be invoked inside a truncated-Newton solver or at the level of processing blocks of 
experimental
data (e.g., per angle). A multigrid approach \cite{lewis2005model} could also
benefit the computational performance.

\section*{Acknowledgments}

We are grateful to Stefan Vogt, Si Chen, Doga Gursoy, and Francesco De
Carlo for introducing us to this problem and for valuable discussions in the
preparation of this paper.
This material was based upon work supported by the U.S. Department of Energy, Office of Science, under contract DE-AC02-06CH11357.

\appendix
\section{Properties of Radon Transform} 
\label{sec:properties_radon} 
We briefly describe a few properties of the Radon transform that explain the equivalence of reconstruction via translation and rotation.\\

{\bf Symmetry:} $Rf(\tau,\theta)=Rf(-\tau,\theta+\pi)$.

{\bf Translation:} Let $g(x-\Delta x,y-\Delta y)=f(x,y)$. Then
\[Rf(\tau,\theta)=Rg(\tau',\theta),\] where $\tau'=\tau-\Delta
x\cos\theta-\Delta y\sin\theta$.  \begin{equation*} \begin{array}{rl}
  Rf(\tau,\theta)&=\ds \int_{-\infty}^\infty \int_{-\infty}^\infty
  f(x,y)\delta(\tau-x\cos\theta-y\sin\theta)dx dy\\ &=\ds
  \int_{-\infty}^\infty \int_{-\infty}^\infty g(x-\Delta x,y-\Delta
  y)\delta(\tau-x\cos\theta-y\sin\theta)dx dy\\ &=\ds \int_{-\infty}^\infty
  \int_{-\infty}^\infty g(x',y')\delta(\tau-(x'+\Delta
  x)\cos\theta-(y'+\Delta y)\sin\theta)dx'dy'\\ &=\ds \int_{-\infty}^\infty
  \int_{-\infty}^\infty
  g(x',y')\delta(\tau'-x'\cos\theta-y'\sin\theta)dx'dy'\\ \vspace{0.1cm} &
  =\ds Rg(\tau',\theta).  \end{array} \end{equation*}

{\bf Rotation:}

Let $g(r,\phi-\Delta\phi)=f(r,\phi)$ where $(r,\phi)$ is the polar
coordinates of $(x,y)$. Then,
\[Rf(\tau,\theta)=Rg(\tau,\theta-\Delta\phi).\] \begin{equation*}
  \begin{array}{rl} Rf(\tau,\theta)& =\ds \int_{-\infty}^\infty
    \int_{-\infty}^\infty
    f(r,\phi)\delta(\tau-r\cos\phi\cos\theta-r\sin\phi\sin\theta)|r|dr
    d\phi\\ &=\ds \int_{-\infty}^\infty \int_{-\infty}^\infty
    f(r,\phi)\delta(\tau-r\cos(\phi-\theta))|r|dr d\phi\\ &=\ds
    \int_{-\infty}^\infty \int_{-\infty}^\infty
    g(r,\phi-\Delta\phi)\delta(\tau-r\cos(\phi-\theta))|r|dr d\phi \\ &=\ds
    \int_{-\infty}^\infty \int_{-\infty}^\infty
    g(r,\phi')\delta(\tau-r\cos(\phi'+\Delta\phi-\theta))|r|dr d\phi' \\
    &=\ds Rg(\tau,\theta-\Delta\phi).  \end{array} \end{equation*}

\section{Derivation of the CoR-Shifted Radon Transform}
\label{SEC:ShiftedRTDerivation}
Here, we present a brief derivation of \eqref{EQN:ShiftedRT}.  Assume, without
loss of generality, that the beam apparatus is initially positioned to measure
the Radon transform along the direction corresponding to $\theta = 0$.  By \Fig{fig:contRadon}, this amounts to integrating the function $f$
representing our object along the vertical lines $x = \tau$, $\tau \in \R$.  To
measure along the direction corresponding to a general angle $\theta$, we
rotate the beam apparatus counterclockwise by $\theta$.%
\footnote{As described earlier, in the
actual experimental setup we rotate the object, not the beam; from a
mathematical point of view, it does not matter which we rotate.}
If the center of rotation is $(x_\theta^*, y_\theta^*)$, the line (beam) $x = \tau$ rotates
to the line
\[
\ell(x_\theta^*, y_\theta^*) = \left\{(\tau - x_\theta^*)\begin{bmatrix}
\cos \theta \\ \sin \theta
\end{bmatrix} + \begin{bmatrix}
x_\theta^* \\ y_\theta^*
\end{bmatrix} + s \begin{bmatrix}
-\sin \theta \\ \cos \theta
\end{bmatrix} \: : \: s \in \R\right\}.
\]
By definition, $Rf(\tau, \theta, x_\theta^*, y_\theta^*)$ is the integral of $f$ along
$\ell(x_\theta^*, y_\theta^*)$.  An easy calculation shows that, ignoring sign, this line is
a distance
\[
\tau' = \tau - x_\theta^* + x_\theta^* \cos \theta + y_\theta^* \sin \theta
\]
from the origin. From the geometry of the rotation process we see that
the line perpendicular to $\ell(x_\theta^*, y_\theta^*)$ is inclined at an angle $\theta$
above the (positive) horizontal axis.  Hence,
\[
Rf(\tau, \theta, x_\theta^*, y_\theta^*) = \int_{\ell(x_\theta^*, y_\theta^*)} f \: |ds| = \int_{-\infty}^\infty \int_{-\infty}^\infty f(x, y) \delta(\tau' - x \cos\theta - y \sin \theta) \: dx \: dy,
\]
from which \eqref{EQN:ShiftedRT} follows immediately upon substituting in for
$\tau'$.

\section{Error Analysis for the Projection Alignment Process}
\label{SEC:ErrorAnalysis}
In this section, we consider the error in the approximation $Rf(\tau, \theta,
0, 0) \approx \widetilde{Rf}(\tau, \theta, 0, 0)$ described in
\Sec{SSEC:Aligning}.  Recall that $\widetilde{Rf}(\tau, \theta, 0, 0)$
incorporates both the Gaussian regularization of the Dirac delta and the
discretization of the convolution integral using the discrete Fourier
transform.

Moving to a more abstract setting, the question is fundamentally one of
assessing the error in the approximation $f(x - P) \approx (f \ast \delta_{P,
\sigma})(x)$, where $f$ is some compactly supported function 
on $\R$ 
and where 
we evaluate the right-hand side by discretizing the convolution integral in a
particular way.  We lose no generality in assuming that $P = 0$, so we can
actually analyze $f(x) \approx (f \ast \delta_\sigma)(x)$, where we use the
shorthand $\delta_\sigma = \delta_{0, \sigma}$. Since we wish to measure the
error with respect to the uniform norm, we need only to consider the
error at a single point $x$.  We lose no generality in assuming that point is
$x = 0$, so we are led to the problem of understanding the error in
\begin{equation}\label{EQN:EAIntegral}
f(0) \approx \int_{-\infty}^\infty f(x) \delta_\sigma(x) \: dx,
\end{equation}
again, under a specific discretization of the integral on the right-hand side,
which in particular includes truncating the domain of the Gaussian
$\delta_\sigma$.

Thus, we assume that $\supp f \subset [-L, L]$ for some $L > 0$, and we define
\[
\delta_\sigma^\mathrm{T}(x ; A) = \begin{cases}
\delta_\sigma(x) & |x| \leq A \\
0                & |x| > A,
\end{cases}
\]
which is $\delta_\sigma$ with its support truncated to $[-A, A]$.  Our goal is
to understand the error $E(\sigma, A, K)$ in the approximation to $f(0)$
obtained by using the $K$-point (composite) trapezoid rule to discretize the
integral in \eqref{EQN:EAIntegral} on the interval $[-M, M]$, where $M = L +
A$:
\[
E(\sigma, A, K) = f(0) - \frac{2M}{K} \sum_{k = -N}^N f(x_k)\delta_\sigma^\textrm{T}(x_k ; A),
\]
where $K = 2N + 1$ is odd%
\footnote{We make this choice solely for notational convenience.  Everything we
do applies for even $K$ as well.}
and $x_k = k(2M/K)$, $-N \leq k \leq N$.  (This is mathematically equivalent to
evaluating the convolution at one of the points $x_k$ using the discrete
Fourier transform.)  In particular, we want to know how $E(\sigma, A, K)$
depends on $\sigma$ and $K$.

$E(\sigma, A, K)$ comprises 
three sources of error.
The first is
the \emph{regularization error} that we incur because $\sigma \neq 0$:
\[
E_\textrm{R}(\sigma) = f(0) - \int_{-\infty}^\infty f(x) \delta_\sigma(x) \: dx.
\]
The second is the \emph{truncation error} that comes from our limiting the
support of the Gaussian:
\[
E_\textrm{T}(\sigma, A) = \int_{-\infty}^\infty f(x) \delta_\sigma(x) \: dx - 
\int_{-\infty}^\infty f(x) \delta_\sigma^\textrm{T}(x ; A) \: dx.
\]
The third is the \emph{discretization error} due to our quadrature rule:
\[
E_\textrm{D}(\sigma, A, K) = \int_{-\infty}^\infty f(x) \delta_\sigma^\textrm{T}(x ; A) \: dx - \frac{2M}{K} \sum_{k = -N}^N f(x_k)\delta_\sigma^\textrm{T}(x_k ; A).
\]
We have
\begin{equation}\label{EQN:CombiningTheErrors}
E(\sigma, A, K) = E_\textrm{R}(\sigma) + E_\textrm{T}(\sigma, A) + E_\textrm{D}(\sigma, A, K).
\end{equation}
We will bound $E(\sigma, A, K)$ by bounding each of the terms on the right-hand
side of \eqref{EQN:CombiningTheErrors} in turn.

Assume that $f$ is twice continuously differentiable.  Then, 
Taylor's theorem (with Lagrange's form of
the remainder) enables us to write
\[
f(x) = f(0) + f'(0) x + \frac{1}{2}f''\bigl(\xi(x)\bigr) x^2
\]
for $x \in \R$, where $\xi(x)$ is some function that satisfies $0 \leq \xi(x)
\leq x$ for $x > 0$ and $x \leq \xi(x) \leq 0$ for $x < 0$.  
From symmetry of $x\delta_\sigma(x)$, it thus follows that
\[
\int_{-\infty}^\infty f(x) \delta_\sigma(x) \: dx = f(0) + \frac{1}{2} \int_{-\infty}^\infty f''\bigl(\xi(x)\bigr) x^2 \delta_\sigma(x) \: dx,
\]
and hence we have
\begin{equation}\label{EQN:RErrorBound}
|E_\textrm{R}(\sigma)| \leq \frac{1}{2} \sigma^2 \|f''\|_\infty
\end{equation}
for the regularization error.

The truncation error is even easier to handle.  We have
\[
E_\textrm{T}(\sigma, A) = \int_{-\infty}^{-A} f(x) \delta_\sigma(x) \: dx + 
\int_A^\infty f(x) \delta_\sigma(x) \: dx,
\]
and so
\begin{equation}\label{EQN:TErrorBound}
|E_\textrm{T}(\sigma, A)| \leq 2 \|f\|_\infty \int_A^\infty \delta_\sigma(x) \: dx \leq \|f\|_\infty e^{-\frac{A^2}{2\sigma^2}},
\end{equation}
by one of the many available bounds for the complementary error function
\cite{chang2011chernoff}.

The discretization error is trickier.  Because $\delta_\sigma^\textrm{T}$ is
discontinuous on $[-M, M]$, we cannot apply most standard results, since these 
results
assume continuity of the integrand.  Instead, we appeal to the following
result from \cite{dragomir2001trapezoid}, which bounds the error in the
trapezoid rule approximation to the integral of a function of bounded
variation.

\begin{proposition}[{\cite[Theorem\ 2.1]{dragomir2001trapezoid}}]
If $f : [a, b] \to \R$ is of bounded variation, then
\[
\left| \int_a^b f(t) \: dt - \frac{b - a}{2}\bigl(f(a) + f(b)\bigr) \right| \leq \frac{b - a}{2} V(f)
\]
where $V(f)$ is the total variation of $f$ over $[a, b]$.
\end{proposition}

This result immediately yields a bound of $(h/2)V(f)$ on the error in the
composite trapezoid rule, where $h$ is the spacing between the grid points.

Because $f$ is twice continuously differentiable, it is of bounded variation, 
and so is
$\delta_\sigma^\textrm{T}$.  Hence, their product is of bounded variation as
well. Moreover, we have
\[
V(f\delta_\sigma^\textrm{T}) \leq \|f\|_\infty V(\delta_\sigma^\textrm{T}) + \|\delta_\sigma^\textrm{T}\|_\infty V(f).
\]
(The norms and total variations are taken over $[-M, M]$.)  Since
\[
\|\delta_\sigma^\textrm{T}\|_\infty = \delta_\sigma^\textrm{T}(0 ; A) = \frac{1}{\sigma \sqrt{2\pi}}
\]
and
\[
V(\delta_\sigma^\textrm{T}) = 2\delta_\sigma^\textrm{T}(0 ; A) = \frac{2}{\sigma \sqrt{2\pi}},
\]
it follows that
\[
V(f\delta_\sigma^\textrm{T}) \leq \frac{1}{\sigma \sqrt{2\pi}} \bigl(2\|f\|_\infty + V(f)\bigr).
\]
The grid spacing in our trapezoid rule approximation is $h = 2M/K$.  It
follows that
\begin{equation}\label{EQN:DErrBound}
|E_\textrm{D}(\sigma, A, K)| \leq \frac{L + A}{\sigma K \sqrt{2\pi}} 
\bigl(2\|f\|_\infty + V(f)\bigr).
\end{equation}

Combining \eqref{EQN:CombiningTheErrors}--\eqref{EQN:DErrBound}, we arrive at
the bound
\[
|E(\sigma, A, K)| \leq \frac{1}{2} \sigma^2 \|f''\|_\infty + \|f\|_\infty 
e^{-\frac{A^2}{2\sigma^2}} + \frac{L + A}{\sigma K \sqrt{2\pi}} 
\bigl(2\|f\|_\infty + V(f)\bigr).
\]
This shows that
\[
|E(\sigma, A, K)| = O(\sigma^2) + O\left(\frac{1}{\sigma K}\right),
\]
the implied limits in the big-O symbols being $\sigma \to 0$ and $K \to
\infty$.  This is sufficient to yield the desired qualitative understanding of
the results presented in \Sec{SSEC:Aligning}.

\bibliographystyle{siam.bst}



\end{document}